\documentclass[a4paper]{amsart}
\usepackage{amssymb,amsfonts}
\usepackage[all,arc]{xy}
\usepackage{enumerate}
\usepackage{mathrsfs}
\usepackage{tikz}
\usepackage{tikz-cd}
\usetikzlibrary{arrows,snakes,backgrounds}
\tikzset{>=stealth}
\usepackage{color}   
\usepackage{marginnote}
\usepackage{sseq}
\usepackage{nameref}
\usepackage[backref=page]{hyperref}
\usepackage{MnSymbol} 
\hypersetup{
    colorlinks=true, 
    linktoc=all,     
    citecolor=black,
    filecolor=black,
    linkcolor=blue,
    urlcolor=black
}
\usepackage{cleveref}
\usetikzlibrary{arrows}

\newtheorem{thm}{Theorem}[section]
\newtheorem*{thm*}{Theorem}

\newtheorem{problem}[thm]{Problem}
\newtheorem*{problem*}{Problem}

\theoremstyle{definition}

\theoremstyle{remark}
\newtheorem{rem}[thm]{Remark}

\newcommand{\bC}{\mathbb{C}}

\newcommand{\bF}{\mathbb{F}}

\newcommand{\bQ}{\mathbb{Q}}

\newcommand{\bZ}{\mathbb{Z}}

\newcommand\Diff{\mathrm{Diff}}

\newcommand\BDiff{\mathrm{BDiff}}

\newcommand\dDiff{\mathrm{Diff}^{\delta}}

\newcommand\BdDiff{\mathrm{BDiff}^{\delta}}

\newcommand\Ker{\operatorname*{Ker}}

\newcommand{\hcoker}{/\!\!/}

\addtocontents{toc}{\protect\setcounter{tocdepth}{1}}
\makeatletter
\let\c@equation\c@thm
\makeatother
\numberwithin{equation}{section}
\title{On powers of the Euler class for flat circle bundles}

\author{Sam Nariman}
\email{sam@math.northwestern.edu}
\address{Department of Mathematics\\
Northwestern University\\
Evanston, IL 60208}

\begin{document}

\begin{abstract}
Apparently a lost theorem of Thurston (\cite{slide}) states that the cube of the Euler class $e^3\in H^6(\BdDiff_{\omega}(S^1);\bQ)$ is zero where $\dDiff_{\omega}(S^1)$ is the analytic orientation preserving diffeomorphisms of the circle with the discrete topology. This is in contrast with Morita's theorem (\cite{morita1984nontriviality}) that the powers of the Euler class are nonzero in $H^*(\BdDiff(S^1);\bQ)$ where $\dDiff(S^1)$ is the orientation preserving $C^{\infty}$- diffeomorphisms of the circle with the discrete topology. The purpose of this short note is to prove that the powers of the Euler class $e^k \in H^*(\BdDiff_{\omega}(S^1);\bZ)$ in fact are nonzero in cohomology with integer coefficients. We also give a short proof of Morita's theorem (\cite{morita1984nontriviality}).
\end{abstract}
\maketitle
\section{introduction}

Let $\Diff(S^1)$ denote the orientation preserving $C^{\infty}$-diffeomorphisms of $S^1$. Let $f: M\to \BDiff(S^1)\simeq \bC P^{\infty}$ be the classifying map for a circle bundle $S^1\to E\to M$ over a manifold $M$. If this bundle is flat then it is induced by a group homomorphism 
\[
\rho: \pi_1(M)\to \Diff(S^1),
\]
which is called the holonomy of the flat bundle. This holonomy homomorphism gives a factorization of the classifying map $f$ as follows
\[
M\to \mathrm{B}\pi_1(M)\xrightarrow{\mathrm{B}\rho} \BdDiff(S^1)\xrightarrow{\iota} \BDiff(S^1),
\]
where $\delta$ means the same group equipped with discrete topology and the map $\iota$ is induced by the identity homomorphism.  Note that $ H^*( \BDiff(S^1);\bZ)=\bZ[e]$ where $e$ in dimension $2$ generates a free polynomial algebra. The generator $e$ is known as the universal Euler class.
\begin{thm}[Morita  \cite{morita1984nontriviality}] \label{Morita} The  map 
$$\iota^*: H^*( \BDiff(S^1);\bQ)\to  H^*(\BdDiff(S^1);\bQ)$$
induces an injection in all degrees.
\end{thm}
\begin{rem}
Our short proof of Morita's theorem imply that $\iota^*$ is injective not only on rational cohomology but also on cohomology with any coefficients.
\end{rem} Geometrically this injection means that for every $k>0$, there exists a manifold $M$ and a holonomy map $\rho: \pi_1(M)\to \Diff(S^1)$ so that the $k$-th power of the Euler class of the associated flat circle bundle, $e(\rho)^k\in H^{2k}(M;\bQ)$ is nonzero. Now let $\Diff_{\omega}(S^1)$ denote the subgroup of analytic orientation preserving diffeomorphisms of the circle. 
\begin{problem*}[A lost theorem of Thurston \cite{slide}]
Prove that if the holonomy group $\rho(\pi_1(M))$ of the flat circle bundle $S^1\to E\to M$ lies in $\Diff_{\omega}(S^1)$ then $e(\rho)^k\in H^{2k}(M;\bQ)$ is zero for $k>2$. In other words, consider the map
\[
\eta: \BdDiff_{\omega}(S^1)\to \BDiff(S^1),
\]
which is induced by the composition $\dDiff_{\omega}(S^1)\hookrightarrow \dDiff(S^1)\to \Diff(S^1)$. Prove that $\eta^*(e)^3=0$ in $H^6(\BdDiff_{\omega}(S^1);\bQ)$.
\end{problem*}
\begin{rem}
It is still not known whether $\eta^*(e^2)\in H^4(\BdDiff_{\omega}(S^1);\bQ)$ is nontrivial (see \cite{moritaopenproblems}).
\end{rem}
Inspired by the work of Milnor on the isomorphism conjecture (see \cite[Theorem 1]{milnor1983homology}) we show that:
\begin{thm}\label{thm1}
The induced map
\[
\eta^*: H^*(\BDiff(S^1);\bZ)\to H^*(\BdDiff_{\omega}(S^1);\bZ)
\]
is injective. In particular, for all positive integers $k$, the class $$\eta^*(e)^k \in H^{2k}(\BdDiff_{\omega}(S^1);\bZ)$$ is not zero.
\end{thm}
\begin{rem}
In \cite{slide} inadvertently  the lost theorem of Thurston is formulated as proving  $\eta^*(e)^3=0$ in $H^6(\BdDiff_{\omega}(S^1);\bZ)$ which was the main motivation for the author to write this note.
\end{rem}
\subsection*{Acknowledgment} I would like to thank \'Etienne Ghys, Takashi Tsuboi and Shigeyuki Morita for the correspondences regarding  Thurston's theorem and Benson Farb for the encouragement to write this note. 

\section{Proof of \Cref{thm1}}
Fix a prime $p$. Let $\bF_p$ be the finite field with $p$-elements. If we embed $\bF_p$ into the circle $S^1$ so that its generator acts by $e^{2\pi i/p}$, we obtain a map
\[
\tau: \mathrm{B}\bF_p\to \mathrm{B}S^1\simeq \bC P^{\infty}.
\]
Recall that the powers of the first Chern class of the complex line bundle classified by $\tau$ are the generators of $H^*(\bF_p;\bF_p)$ in the even degrees, therefore the map $\tau^*$ induces an injection
\[
\tau^*: H^*(\mathrm{B}S^1;\bF_p)\hookrightarrow H^*( \mathrm{B}\bF_p;\bF_p).
\]
Now let us think of $S^1$ as a subgroup of rotation matrices in $\Diff_{\omega}(S^1)$. Therefore we have a homotopy commutative diagram 
 
   \[
   \begin{tikzpicture}[node distance=2.8cm, auto]
  \node (A) {$\BdDiff_{\omega}(S^1)$};
  \node (B) [right of=A] {$\BDiff(S^1)$};
  \node (C) [below of=A, node distance=1.7cm] {$\mathrm{B}\bF_p$};  
  \node (D) [below of=B, node distance=1.7cm] {$\mathrm{B}S^1$,};
  \draw[->] (C) to node {$\tau$} (D);
  \draw [->] (A) to node {$\eta$} (B);
    \draw [<-] (A) to node {} (C);
  \draw [<-] (B) to node {$\simeq$} (D);
\end{tikzpicture}
 \]
Since $\tau^*$ is injective in  the following commutative diagram on cohomology

 \[
   \begin{tikzpicture}[node distance=4.3cm, auto]
  \node (A) {$H^*(\BdDiff_{\omega}(S^1);\bF_p)$};
  \node (B) [right of=A] {$H^*(\BDiff(S^1);\bF_p)$};
  \node (C) [below of=A, node distance=1.7cm] {$H^*(\mathrm{B}\bF_p;\bF_p)$};  
  \node (D) [below of=B, node distance=1.7cm] {$H^*(\mathrm{B}S^1;\bF_p)$,};
  \draw[<-] (C) to node {$\tau^*$} (D);
  \draw [<-] (A) to node {$\eta_p^*$} (B);
    \draw [->] (A) to node {} (C);
  \draw [->] (B) to node {$\cong$} (D);
\end{tikzpicture}
 \]
for every prime $p$, the induced map
\[
\eta_p^*:H^*(\BDiff(S^1);\bF_p)\hookrightarrow H^*(\BdDiff_{\omega}(S^1);\bF_p),
\]
is also injective. Now to show that $\eta^*$ is also injective on cohomology with integer coefficients, consider the induced map between Bockstein exact sequences
 \[
   \begin{tikzpicture}[node distance=4.3cm, auto]
  \node (A) {$H^*(\BDiff(S^1);\bF_p)$};
  \node (B) [right of=A] {$H^*(\BDiff(S^1);\bZ)$};
    \node (E) [right of=B] {$H^*(\BDiff(S^1);\bZ)$};
  \node (F) [right of=D] {$H^*(\BdDiff_{\omega}(S^1);\bZ).$};
  \node (C) [below of=A, node distance=1.7cm] {$H^*(\BdDiff_{\omega}(S^1);\bF_p)$};  
  \node (D) [below of=B, node distance=1.7cm] {$H^*(\BdDiff_{\omega}(S^1);\bZ)$};
  \draw[<-] (C) to node {$$} (D);
  \draw [<-] (A) to node {$$} (B);
    \draw [->] (A) to node {$\eta^*_p$} (C);
  \draw [->] (B) to node {$\eta^*$} (D);
    \draw [->] (E) to node {$\eta^*$} (F);
    \draw [<-] (B) to node {$\times p$} (E);
    \draw [<-] (D) to node {$\times p$} (F);
\end{tikzpicture}
 \]
Suppose for some $a\in \Ker(\eta^*)$. Since $\eta^*_p$ is injective for all prime $p$, we have $a\in pH^*(\BDiff(S^1);\bZ)=p\bZ[e]$ for all prime $p$. Therefore $a=0$.\qed
\begin{rem}
The above argument also implies that the map
\[
\iota^*: H^*(\BDiff(S^1);\bZ)\to H^*(\BdDiff(S^1);\bZ)
\]
is injective in all degrees. But note that since $\BdDiff(S^1)$ is not a finite type space, the universal coefficient theorem does not imply \Cref{Morita}. 
\end{rem}
Using different methods, the author proved similar results for diffeomorphisms of other manifolds namely surfaces, higher dimensional analogue of surfaces and a punctured $2$-disk (\cite{nariman2015stable, nariman2014homologicalstability, nariman2015braid}). For instance, for the case of surfaces, let $\Sigma$ denote an orientable surface with genus $g(\Sigma)$ and let $\Diff_c(\Sigma)$ be the group of orientation preserving diffeomorphisms of $\Sigma$ whose supports are in compact subsets of the interior of $\Sigma$ if it has  boundary.  Then the map
\[
\iota^*: H^*(\BDiff_c(\Sigma);\bZ)\to H^*(\BdDiff_c(\Sigma);\bZ)
\]
is injective for $*\leq (2g(\Sigma)-2)/3$. 

 It is interesting to see if there is a general statement for diffeomorphism groups of manifolds as topological groups similar to Milnor's theorem for Lie groups \cite[Corollary 1]{milnor1983homology}.
\begin{problem}
Let $M$ be a manifold with or without boundary. Let $\Diff_c(M)$ denote the group of diffeomorphisms whose supports are away from the boundary of $M$ if it has boundary. Is the induced map
\[
H^*(\BDiff_c(M);\bZ)\to H^*(\BdDiff_c(M);\bZ),
\]
 injective?
\end{problem}
\section{A short proof of Morita's theorem}
Let $\overline{\BDiff(S^1)}$ denote the homotopy fiber of the map $\iota$
\[
\overline{\BDiff(S^1)}\to \BdDiff(S^1)\xrightarrow{\iota}\BDiff(S^1).
\]
To choose a model for the homotopy fiber, let $\mathrm{E}\Diff(S^1)$ denote the universal $\Diff(S^1)$-bundle over $\BDiff(S^1)$. The pullback of $\mathrm{E}\Diff(S^1)$ over $\BdDiff(S^1)$ is a model for  $\overline{\BDiff(S^1)}$. Therefore the group $\Diff(S^1)$ acts on $\overline{\BDiff(S^1)}$ and the homotopy quotient (Borel construction) of this action is weakly homotopy equivalent to
\[
\overline{\BDiff(S^1)}\hcoker \Diff(S^1)\xrightarrow{\simeq} \BdDiff(S^1).
\]
 Since the inclusion $S^1\hookrightarrow \Diff(S^1)$ is a homotopy equivalence, we have the induced maps between fibrations
\[
   \begin{tikzpicture}[node distance=3cm, auto]
  \node (A) {$\overline{\BDiff(S^1)}$};
  \node (B) [below of=A, node distance=1.5cm] {$\overline{\BDiff(S^1)}\hcoker S^1$};
  \node (C) [right of=A, ] {$\overline{\BDiff(S^1)}$};  
  \node (D) [below of=C, node distance=1.5cm] {$\BdDiff(S^1)$};
    \node (E) [below of=B, node distance=1.5cm] {$\mathrm{B}S^1$};
  \node (F) [below of=D, node distance=1.5cm] {$\BDiff(S^1)$};
  \draw[->] (C) to node {$$} (D);
  \draw [->] (A) to node {$$} (B);
    \draw [->] (A) to node {$=$} (C);
  \draw [->] (B) to node {$\simeq$} (D);
    \draw [->] (E) to node {$\simeq$} (F);
  \draw [->] (B) to node {$\pi$} (E);
  \draw [->] (D) to node {$$} (F);
\end{tikzpicture}
 \]
such that every horizontal map induces a weak homotopy equivalence. Therefore, it is enough to show that $\pi$ induces an injection on cohomology with $\bQ$-coefficients.

To do so, we use Thurston's theorem ( \cite[Theorem 5]{thurston1974foliations}) to replace $\overline{\BDiff(S^1)}$ with a free loop space. To recall Thurston's theorem, let $\mathrm{BS}\Gamma_1$ be the classifying space of the orientable codimension $1$ Haefliger structures (see \cite{haefliger1971homotopy} for definitions). Let $\theta$ denote the normal bundle of the universal $\mathrm{S}\Gamma_1$-structure on $\mathrm{BS}\Gamma_1$.  Since it is an orientable line bundle over $\mathrm{BS}\Gamma_1$, it is trivializable. 

 We denote  by $\text{Bun}(TS^1,\theta)$ the space of all bundle maps $TS^1\rightarrow \theta$ from the tangent bundle of $S^1$ to $\theta$  equipped with the compact-open topology.  The action of $\Diff(S^1)$  is given by precomposing a bundle map with the differential of a diffeomorphism.

In  \cite[Section 5.1]{nariman2014homologicalstability}, we reformulated  Thurston's theorem as follows:
\begin{thm}
There is a $\Diff(S^1)$-equivariant map 
\[
f: \overline{\BDiff(S^1)}\to \text{\textnormal{Bun}}(TS^1,\theta)
\]
that induces a homology isomorphism.
\end{thm}

Note that since $TS^1$  is a trivial bundle and the trivialization is equivariant for $\Diff(S^1)$, the space $\text{Bun}(TS^1,\theta)$ is weakly equivalent to $\text{\textnormal{Map}}(S^1,\mathrm{BS}\Gamma_1)$. Similar to \cite[Section 5.1]{nariman2014homologicalstability}, in fact there is a map 
\[
\alpha: \text{\textnormal{Map}}(S^1,\mathrm{BS}\Gamma_1)\to \text{Bun}(TS^1,\theta),
\]
that induces a weak equivalence. Note that the group $S^1$ acts on the both sides of the map $\alpha$. It acts on the free loop space $\text{\textnormal{Map}}(S^1,\mathrm{BS}\Gamma_1)$ by rotating the domain of loops and it acts on $ \text{Bun}(TS^1,\theta)$ as a subgroup of $\Diff(S^1)$. Since the action of $S^1$ on  $\text{Bun}(TS^1,\theta)$ does not change the length of the tangent vectors, one can easily see that the map $\alpha$ is $S^1$-equivariant. Therefore, we have a zig-zag of weak-equivalences
\[
\text{\textnormal{Map}}(S^1,\mathrm{BS}\Gamma_1)\hcoker S^1\to  \text{Bun}(TS^1,\theta)\hcoker S^1\leftarrow \overline{\BDiff(S^1)}\hcoker S^1.
\]
Hence it is enough to show that the map 
\[
p: \text{\textnormal{Map}}(S^1,\mathrm{BS}\Gamma_1)\hcoker S^1\to \mathrm{B}S^1
\]
induces an injection on cohomology with rational coefficients. Since the action of $S^1$ on the free loop space has fixed points (e.g. constant loops), the map $p$ has a section. Thus, the induced map on cohomology with any coefficients is injective and in particular
\[
p^*: H^*( \mathrm{B}S^1;\bQ)\hookrightarrow H^*( \text{\textnormal{Map}}(S^1,\mathrm{BS}\Gamma_1)\hcoker S^1;\bQ)
\]
is injective.\qed
\begin{rem}
Morita informed the author that Haefliger also simplified his argument using  Sullivan's minimal model for the Borel constructions. After all, using rational homotopy theory does not seem to be necessary to prove Morita's theorem.
\end{rem}
\bibliographystyle{abbrv}
\bibliography{reference}
\end{document}